\documentclass{elsart3-1}
\usepackage{amsmath}
\usepackage{amssymb,latexsym}
\usepackage{amscd}
\usepackage{xspace}
\usepackage{verbatim}

\usepackage[english,francais]{babel}

\newtheorem{e-proposition}[theorem]{Proposition}

\newtheorem{e-definition}[theorem]{Definition\rm}

\renewcommand{\theequation}{\arabic{equation}}
\def\og{\leavevmode\raise.3ex\hbox{$\scriptscriptstyle\langle\!\langle$~}}
\def\fg{\leavevmode\raise.3ex\hbox{~$\!\scriptscriptstyle\,\rangle\!\rangle$}}
\newcommand{\txt}[1]{\;\text{ #1 }\;}

\newcommand{\be}{\begin{equation}}
\newcommand{\bel}[1]{\begin{equation}\label{#1}}
\newcommand{\ee}{\end{equation}}
\newtheorem{subn}{\name}

\newcommand{\bsn}[1]{\def\name{#1$\!\!$}\begin{subn}}
\newcommand{\esn}{\end{subn}}
\newtheorem{sub}{\name}[section]
\newcommand{\bs}{\begin{sub}}
\newcommand{\es}{\end{sub}}
\newcommand{\bsl}[1]{\begin{sub}\label{#1}}
\newcommand{\bth}[1]{\def\name{Theorem}\begin{sub}\label{t:#1}}
\newcommand{\blemma}[1]{\def\name{Lemma}\begin{sub}\label{l:#1}}
\newcommand{\bcor}[1]{\def\name{Corollary}\begin{sub}\label{c:#1}}
\newcommand{\bdef}[1]{\def\name{Definition}\begin{sub}\label{d:#1}}
\newcommand{\bprop}[1]{\def\name{Proposition}\begin{sub}\label{p:#1}}
\newcommand{\bnote}[1]{\def\name{\mdseries\scshape Notation}\begin{sub}\label{n:#1}}
\newcommand{\bcom}{}
\newcommand{\req}{\eqref}
\newcommand{\rth}[1]{Theorem~\ref{t:#1}}
\newcommand{\rlemma}[1]{Lemma~\ref{l:#1}}

\newcommand{\BA}{\begin{array}}
\newcommand{\EA}{\end{array}}
\newcommand{\BAN}{\renewcommand{\arraystretch}{1.2}
\setlength{\arraycolsep}{2pt}\begin{array}}
\newcommand{\BAV}[2]{\renewcommand{\arraystretch}{#1}
\setlength{\arraycolsep}{#2}\begin{array}}
\newcommand{\BSA}{\begin{subarray}}
\newcommand{\ESA}{\end{subarray}}
\newcommand{\BAL}{\begin{aligned}}
\newcommand{\EAL}{\end{aligned}}
\newcommand{\BALG}{\begin{alignat}}
\newcommand{\EALG}{\end{alignat}}
\newcommand{\BALGN}{\begin{alignat*}}
\newcommand{\EALGN}{\end{alignat*}}

\newcommand{\forevery}{\quad \forall}
\newcommand{\1}{\\[1mm]}
\newcommand{\2}{\\[2mm]}

\newcommand{\set}[1]{\{#1\}}
\def\({{\rm (}}
\def\){{\rm )}}

\def\uar{\uparrow}
\newcommand{\abs}[1]{\left |#1\right |}

\newcommand{\opname}[1]{\mathrm{#1}\,}

\newcommand{\dist}{\opname{dist}}

\newcommand{\q}{\quad}

\newcommand{\prt}{\partial}
\newcommand{\sms}{\setminus}

\newcommand{\ti}{\times}

\newcommand{\tl}{\tilde}

\newcommand{\sbs}{\subset}

\newcommand{\tin}{\to\infty}
\newcommand{\ind}[1]{_{_{#1}}\!}
\newcommand{\chr}[1]{\chi\ind{#1}}

\newcommand{\num}[1]{{\rm (#1)}\hspace{2mm}}

\newcommand{\wrto}{with respect to\xspace}

\newcommand{\consy}{consequently\xspace}

\newcommand{\ngh}{neighborhood\xspace}

\newcommand{\seq}{sequence\xspace}

\newcommand{\ifif}{if and only if\xspace}

\newcommand{\abc}{absolutely continuous\xspace}

\newcommand{\bdw}{\partial\Gw}
\newcommand{\Capq}{C_{2/q,q'}}

\newcommand{\sbsq}{\overset{q}{\sbs}}

\newcommand{\tr}{\mathrm{tr}\,}

\newcommand{\qcl}{$q$-closed\xspace}
\newcommand{\qop}{$q$-open\xspace}
\newcommand{\gsmod}{$\gs$-moderate\xspace}

\newcommand{\ugb}[1]{u\chr{\Gs_\gb(#1)}}
\newcommand{\mcon}{$q$-moderately convergent\xspace}
\newcommand{\mdiv}{$q$-moderately divergent\xspace}

\def\bcom{}
\def\ga{\alpha}     \def\gb{\beta}

\def\gm{\mu}        \def\gn{\nu}         
    \def\gr{\rho}        
\def\gs{\sigma}       
      \def\gw{\omega}
\def\gx{\xi}                
     \def\Gd{\Delta}      

    \def\Gs{\Sigma}      
\def\Gw{\Omega}              

\def\CS{{\mathcal S}}      
\def\CR{{\mathcal R}}      
      
      \def\CF{{\mathcal F}}

   \def\CU{{\mathcal U}}

   \def\BBR {\mathbb R}

\def\GTT {\mathfrak T}

\def\rqq{\req{q-eq}\xspace}

\begin{document}
\centerline{}
\begin{frontmatter}


\selectlanguage{english}
\title{The precise boundary trace of solutions of
a class of supercritical nonlinear equations}


\selectlanguage{english}
\author[authorlabel1]{Moshe Marcus},
\ead{marcusm@math.technion.ac.il}
\author[authorlabel2]{Laurent V\'eron}
\ead{veronl@univ-tours.fr}

\address[authorlabel1]{Department of Mathematics, Technion\\
 Haifa 32000, ISRAEL}
\address[authorlabel2]{Laboratoire de Math\'ematiques et Physique Th\'eorique, Facult\'e des Sciences\\
Parc de Grandmont, 37200 Tours, FRANCE}


\medskip

\begin{abstract}
\selectlanguage{english}
We construct and study the properties of the precise boundary trace of positive solutions of $-\Delta u+u^q=0$ in a smooth bounded domain of $\mathbb R^N$, in the supercritical case $q\geq q_c=(N+1)/(N-1)$.
{\it To cite this article: 
M. Marcus, L. V\'eron, C. R. Acad. Sci. Paris, Ser. I XXX (2006).}

\vskip 0.5\baselineskip

\selectlanguage{francais}
\noindent{\bf R\'esum\'e} \vskip 0.5\baselineskip \noindent
{\bf La trace au bord pr\'ecise des solutions d'une classe d'\'equations non lin\'eaires sur-critiques}
Nous construisons et \'etudions les propri\'et\'es de la trace au bord pr\'ecise des solutions positives de $-\Delta u+u^q=0$ dans un domaine 
r\'egulier de $\mathbb R^N$, dans le cas sur-critique $q\geq q_c=(N+1)/(N-1)$. {\it Pour citer cet article~: M. Marcus, L. V\'eron, C. R. Acad. Sci. Paris, Ser. I XXX (2006).}
\end{abstract}
\end{frontmatter}
\selectlanguage{francais}
\section*{Version fran\c{c}aise abr\'eg\'ee}
\setcounter {equation}{0}
Soit $\Gw$ un ouvert born\'e de $\BBR^N$ de bord de classe $C^2$ et $u\in L^q_{loc}(\Gw)$ ($q>1$) une solution positive de 
\begin{equation}\label{Fq-eq}
 -\Gd u+|u|^{q-1}u=0.
\end{equation}
Il est bien connu que $u$ poss\`ede une trace au bord dans la classe des mesures de Borel ayant la r\'egularit\'e ext\'erieure au sens classique. Si $q\geq q_{c}$ cette notion de trace n'est pas suffisante pour d\'eterminer de fa\c{c}on unique la solution de (\ref{Fq-eq}). Contrairement \`a la trace fine \cite {DK3}, qui s'exprime en termes probabilistes et est limit\'ee au cas $q_{c}\leq q\leq 2$, la notion de trace pr\'ecise que  nous d\'eveloppons est valable pour tout $q\geq q_{c}$. Elle est fond\'ee sur la topologie fine $\frak T_{q}$ associ\'ee \`a la capacit\'e de Bessel $C_{{2/q,q'}}$ sur $\prt\Gw$. Notons $\gr(x):=\dist (x,\prt\Gw)$ et 
$\Gw_{\gb}=\{x\in\Gw:\gr(x)<\gb\},\quad\Gw'_{\gb}=\gw\setminus\bar\Gw_{\gb},\quad\Gs_{\gb}=\prt\Gw'_{\gb}$.
 Il existe $\gb_{0}>0$ tel que pour tout $x\in \Gw_{\gb_{0}}$ il existe un unique $\gs(x)\in\Gw$ tel que 
$\gr(x)=\abs{x-\gs(x)}$. Si $Q\subset\prt\Gw$ est $\frak T_{q}$-ouvert, on note $\Gs_{\gb}(Q)=\{x\in\Gs_{\gb}:\gs(x)\in Q\}$, et, si $u\in C(\Gw)$, $u^Q_{\gb}$ d\'esigne la solution de (\ref{Fq-eq}) dans $\Gw_{\gb}'$ valant $u\chi_{_{\Gs_{\gb}(Q)}}$ sur $\Gs_{\gb}$. La dichotomie suivante est \`a la base de nos r\'esultats.\medskip

\noindent {\bf Th\'eor\`eme 1 }{\it
Soit $u$ une solution positive de (\ref {Fq-eq}). Si $Q\subset\prt\Gw$ est $\frak T_{q}$-ouvert et si la limite suivante existe, on pose $L(Q)=\lim_{\gb\to 0}\int_{\Gs_{\gb}(Q)}udS$. Alors,

 \noindent {\bf   ou bien,} $L(Q)=\infty$ pour tout  voisinage $\GTT_q$-ouvert $Q$ de $\,\gx$,
 
  \noindent {\bf  ou bien,} il existe  un voisinage $\GTT_q$-ouvert $Q$ de $\xi$ tel que $L(Q)<\infty$. 
  
  \noindent Le premier cas se produit si et seulement si, pour tout $\frak T_{q}$- voisinage $Q$ de $\xi$,
  \begin{equation}\label{BU}
  \int_{A}u^q\gr(x)dx=\infty,\quad A=(0,\gb_{0})\ti Q.
  \end {equation}}

Un point $\xi$ est dit {\it singulier} (resp. {\it r\'egulier}) si le premier (resp. le deuxi\`eme) cas se produit. L'ensemble des point singuliers (resp. r\'eguliers) not\'e $\CS(u)$ (resp. $\CR(u)$ est $\GTT_q$-ferm\'e (resp. $\GTT_q$-ouvert). Si $A\subset \prt\Gw$, nous notons $\tilde A$ la fermeture de $A$ dans topologie $\frak T_{q}$.\medskip

\noindent {\bf Th\'eor\`eme 2 }{\it Il existe une mesure de Borel positive $\mu$ sur $\bdw$
poss\'edant les propri\'et\'es suivantes:\1
\num{i} Pour tout $\gs\in
\CR(u)$ il existe un voisinage $\GTT_q$-ouvert $Q$ de $\gs$ et
une solution mod\'er\'ee $w$ de (\ref{Fq-eq}) tels que $\tl Q\sbs \CR(u)$, $\mu(\tl Q)<\infty$ et
 \begin{equation}\label{intlim}
 u^{Q}_\gb\to w \txt{localement uniform\'ement dans $\Gw$}\!\!\!, \quad (\tr w)\chr{Q}= \mu\chr{Q}.
\end{equation}
\num{ii} $\mu$ a la r\'egularit\'e ext\'erieure pour la topologie $\GTT_q$ et est absolument continue par rapport \`a la capacit\'e $C_{2/q,q'}$ sur les sous-ensembles $\GTT_q$-ouverts o\`u elle est born\'ee.
}\medskip

Le couple $(\gm,\CS(u)$ est, par d\'efinition, la {\it trace pr\'ecise} de $u$, not\'ee $tr(u)$, qui peut \^etre aussi repr\'esent\'ee par la mesure de Borel $\gn$ d\'efinie par $\gn=\gm$ sur $\CR(u)$ et $\gn(A)=\infty$ pour tout bor\'elien $A$ tel que $A\cap\CS(u)\neq\emptyset$; $\gn$ a les propri\'et\'es suivantes: (i) Elle a la r\'egularit\'e ext\'erieure pour la topologie $\GTT_q$. (ii) Elle est absolument continue par rapport \`a la capacit\'e $C_{2/q,q'}$ au sens o\`u pour tout ensemble $\GTT_q$-ouvert $Q$ et tout bor\'elien $A$, $C_{2/q,q'}(A)=0$ implique $\gn(Q)=\gn(Q\setminus A)$. Une mesure de Borel v\'erifiant (i) et (ii) est dite {\it $q$-parfaite}. Nous donnons alors la condition n\'ecessaire et suffisante d`existence,
ainsi qu'un r\'esultat d'unicit\'e, de la solution du probl\`eme aux limites g\'en\'eralis\'e
 \begin{equation}\label{pblim}
-\Gd u+u^q=0, \quad u>0\text { dans }\Gw,\quad tr(u)=\gn,
\end{equation}

\noindent {\bf Th\'eor\`eme 3 }{\it Soit $\nu$ une mesure de Borel sur $\bdw$, born\'ee ou non.
 Le probl\`eme aux limites (\ref{pblim})
a une solution si et seulement si  $\nu$ est $q$-parfaite. Quand c'est le cas,
une solution de (\ref{pblim}) est donn\'ee par
\begin{equation}\label{precisely}
 U=v\oplus U_{F}, \q v=\sup\set{u_{\nu\chr{Q}}:Q\in \CF_\nu},
\end{equation}
o\`u
$\CF_\nu:=\set{Q : \, Q\;q\text{-ouvert}, \;\nu(Q)<\infty}$, $G:=\bigcup_{\CF_\nu}Q$,
$F=\bdw\sms G$, $U_F$ est la solution maximale s'annullant sur $\bdw\sms F$ et $v\oplus U_{F}$ est la plus grande solution de (\ref{Fq-eq}) inf\'erieure \`a la sur-solution $v+U_{F}$. Enfin $U$ est $\gs$-mod\'er\'ee, c'est la solution maximale du probl\`eme (\ref{pblim}) dont c'est l'unique solution $\gs$-mod\'er\'ee. }


\selectlanguage{english}
\section{Introduction and statement of results} In this paper we present a theory of boundary trace of positive solutions of the equation
\begin{equation}\label{q-eq}
 -\Gd u+|u|^{q-1}u=0
\end{equation}
 in a  bounded domain $\Gw\sbs \BBR^N$ of class $C^2$. A function $u$ is a solution if
 $u\in L^q_{loc}(\Gw)$ and the equation holds in the distribution sense.
\par Equations of this type
 have been intensively studied in the last ten years in the context of  boundary trace theory
and the associated boundary value problem.
In the subcritical case, $1<q<q_c=(N+1)/(N-1)$, the problem is well understood thanks to  the
works of Le Gall \cite{LG},
using a probabilistic approach which imposes $q\leq 2$, and
Marcus and V\'eron \cite{MV1} by using an analytic approach, with no restriction.
In the supercritical case $q\geq q_{c}$, the notions of removable sets and admissible
Radon measures
 are implemented by Le Gall,  by Dynkin and Kuznetsov \cite{DK1} and by
 Marcus and V\'eron \cite{MV2,MV3}.  But in 1997 Le Gall showed that the standard trace theory is
 not appropriate because many solutions may have the same boundary trace. Following this
 observation, a  theory of
{\it fine trace}, based upon a probabilistic formulation, was introduced by Dynkin and Kuznetsov  \cite{DK3} who showed that, for $q\leq 2$,
 the fine trace theory is satisfactory in the family of so-called \gsmod solutions.
 Later on it was shown by Mselati
\cite{Ms}, combining Le Gall's Brownian snake method \cite{LGbook} and Dynkin's approach
\cite {Dbook1} that, in the case $q=2$, all positive solutions are \gsmod.
By analytical methods Marcus and
V\'eron \cite{MV4} proved that, for all $q\geq q_c$ and every compact set $K\sbs \bdw$,
 the maximal solution of \rqq vanishing on $\bdw \sms K$ is \gsmod. Their
proof
was based on the derivation of sharp capacitary estimates for the
maximal solution. In continuation, Dynkin \cite{Dbook2}, using the estimates of Marcus and
Veron  \cite{MV4}, extended Mselati's result to all case $q\leq 2$.
For $q>2$ the problem remains open.

\par Our definition of boundary trace is based on the fine topology
associated
 with the Bessel capacity $\Capq$ on $\bdw$, denoted by $\GTT_q$. The $\GTT_q$-closure of a set $E$
 will be denoted by $\tl E$. We also need the following notation.\smallskip

\noindent {\it Notation 1.1}
\num{a} For every $x\in \BBR^N$ and every $\gb>0$ put
  $\gr(x):=\dist(x,\prt\Gw)$ and
$\Gw_\gb=\set{x\in \Gw:\,\gr(x)<\gb}$, $\Gw'_\gb=\Gw\sms \bar \Gw_\gb$, $\Gs_\gb=\prt\Gw'_\gb.$\\
\num{b} There exists a positive number $\gb_0$ such that, for every $x\in \Gw_{\gb_0}$ there exists
  a unique point $\gx\in\bdw$ such that $\dist(x,\gx)=\gr(x)$. Put $\gs(x):=\gx$.\\
\num{c} If $Q$ is a $\GTT_q$-open subset of $\bdw$ and $u\in C(\bdw)$ we denote by $u_\gb^Q$
the solution of \rqq in $\Gw'_\gb$ with boundary data
$h=u\chr{\Gs_\gb(Q)}$ where $\Gs_\gb(Q)=\{x\in \Gs_\gb:\, \gs(x)\in Q\}$.\2
\indent Recall
that a  solution $u$ is moderate if $\abs{u}$ is dominated by a
harmonic function. When this is the case, $u$ possesses a boundary
trace  (denoted by $\tr u$) given by a bounded Borel measure.
 A positive solution $u$ is \gsmod if there exists an increasing \seq of moderate solutions $\set{u_n}$
such that $u_n\uar u$. This notion was introduced by Dynkin and
Kuznetsov \cite{DK3} (see also \cite{Ku} and \cite{Dbook1}).

\par If $\nu$ is a bounded Borel measure on $\bdw$, the  problem
\begin{equation}\label{BVP}
 -\Gd u+u^q=0 \text{ in }\Gw, \q \tr u=\nu \text{ on }\bdw
\end{equation}
possesses a (unique) solution \ifif $\nu$ vanishes on sets of
$\Capq$-capacity zero, (see \cite {MV3} and the references
therein). The solution is denoted by $u_\nu$.
\par
 The set of positive solutions of \rqq in $\Gw$ will be denoted by $\CU(\Gw)$. It is well
 known that this set is compact in the topology of $C(\Gw)$, i.e.,
 relative to local uniform convergence in $\Gw$. If $u,v\in \CU(\Gw)$, we denote by $u\oplus v$
 the largest solution dominated by $u+v$.

 The first result displays a dichotomy which is the basis for our definition of boundary trace.\medskip
\bth{dichotomy}
 Let $u\in \CU(\Gw)$  and let $\gx\in \bdw$. If $\,Q\sbs \bdw$ is a $\GTT_q$-open set
and  provided the following limit exists, put $L_Q=\lim_{\gb\to 0}\int_{\Gs_\gb(Q)}udS$.
Then,
\begin{quotation}
{\bf either} $L_Q=\infty$ for every $\GTT_q$-open \ngh $Q$ of $\,\gx$,\smallskip

\noindent{\bf or} there exists  a $\GTT_q$-open \ngh $Q$ such that $L_Q<\infty$.
\end{quotation}
\bcom
\begin{quotation}
  {\bf either,} for every $\GTT_q$-\ngh $Q$ of $\,\gx$, $\lim_{\gb\to 0}\int_{\Gs_\gb(Q)}udS=\infty$\\
{\bf or} there exists  a $\GTT_q$-\ngh $Q$ such that the following limit exists
$ \lim_{\gb\to 0}\int_{\Gs_\gb(Q)}udS<\infty.$
\end{quotation}

\begin{description}
  \item[either] for every $\GTT_q$-open \ngh $Q$ of $\,\gx$,
  $\set{\ugb{Q}}$ is \mdiv as $\gb\to 0$,
\item[or] there exists  a $\GTT_q$-open \ngh $Q$ of $\,\gx$ such that
$\set{\ugb{Q}}$ is \mcon as $\gb\to0$.
\end{description}
\end{comment}
The first case  occurs \ifif, for every $\GTT_q$-\ngh $Q$ of $\gx$,
\begin{equation}\label{int-sing}
\int_{A}u^q\gr(x)dx=\infty,\q A=(0,\gb_0)\ti Q.
\end{equation}
\es

A point $\gx\in\bdw$ is called a {\em singular} point of $u$  in the first case and
 a {\em regular} point otherwise.
The set of singular points is denoted by $\CS(u)$ and its
complement in $\bdw$ by $\CR(u)$.
Our next result provides additional information on the behavior of solutions near the regular
boundary set.\medskip

\bth{int-B} Let $u\in \CU(\Gw)$. Then $\CR(u)$ is $\GTT_q$-open
and there exists a non-negative Borel measure $\mu$ on $\bdw$
possessing the following properties.\1
\num{i} For every $\gs\in
\CR(u)$ there exist a $\GTT_q$-open \ngh $Q$ of $\gs$ and a
moderate solution $w$ such that $\tl Q\sbs \CR(u)$, $\mu(\tl Q)<\infty$ and
 \begin{equation}\label{intlim}
 u^{Q}_\gb\to w \txt{locally uniformly in $\Gw$}, \q (\tr w)\chr{Q}= \mu\chr{Q}.
\end{equation}
\num{ii} $\mu$ is outer regular relative to $\GTT_q$ and  \abc
relative to $\Capq$ on $\GTT_q$-open sets on which it is bounded.
\es

Based on these results we define the {\em precise boundary
trace} of $u$ as the couple $(\mu,\CS(u)).$
The trace can also be represented by a Borel measure
$ \nu$ defined as follows. For every Borel set $A\sbs \bdw$:
\begin{equation}\label{PT}
\nu(A)=\mu(A)\; \text{ if }A\sbs \CR(u),\q \nu(A)=\infty\;\text{ otherwise}.
\end{equation}
We denote $\tr^c=(\mu,\CS(u))$ and $\tr u=\nu$.
\bcom
Let $\nu$ be a positive Borel measure on $\bdw$,
possibly unbounded. We shall say that $\gx$ is a singular point of $\nu$ if $\nu(Q)=\infty$ for every
$\GTT_q$-open \ngh $Q$ of $\,\gx$. The set of singular points of $\nu$ is denote by $\CS_\nu$ and its
complement in $\bdw$ by $\CR_\nu$.
We say that $u\in \CU(\Gw)$ is a solution of \req{BVP} if $\CS(u)=\CS_\nu$ and for every $\gx\in \CR_\nu$ there
exists a $\GTT_q$-open \ngh $Q$ such that $u^Q=\lim_{\gb\to 0}u_\gb^Q$ exists and
is a moderate solution with $\tr u^Q=\nu\chr(Q)$.
\end{comment}
\par Let $\nu$ be a positive Borel measure on $\bdw$. We say that $\mu$ is {\em $q$-perfect} if:\1
(i) $\mu$ is outer regular relative to
$\GTT_q$. (ii) $\mu$ is {\em essentially absolutely continuous}
relative to $\Capq$, i.e., if $Q$ is  $\GTT_q$-open   and  $\Capq(A)=0$ then
$\nu(Q)=\nu(Q\sms A).$\1
The second property  implies that, if $\nu(Q\sms A)<\infty$ then
$\nu(A\cap Q)=0$.

\par We have the following existence and uniqueness results for the
(generalized) boundary value problem \req{BVP}, where $\tr u=\nu$ is understood as in \req{PT}.\medskip

\bth{exist} Let $\nu$ be a positive Borel measure, possibly unbounded. Then
\req{BVP}
possesses a solution \ifif $\nu$ is $q$-perfect. When this
condition holds, a solution of \req{BVP} is given by
\begin{equation}\label{precisely}
 U=v\oplus U_{F}, \q v=\sup\set{u_{\nu\chr{Q}}:Q\in \CF_\nu},
\end{equation}
where
$\CF_\nu:=\set{Q : \, Q\textrm{ \qop }, \;\nu(Q)<\infty}$, $G:=\bigcup_{\CF_\nu}Q$, $F=\bdw\sms G$
and $U_F$ is the maximal solution vanishing on $\bdw\sms F$. \es
\medskip

\bth{unique}
Let $\nu$ be a  $q$-perfect measure on $\bdw$. Then the solution
$U$ of problem \req{BVP} defined by \req{precisely} is \gsmod
and it is the maximal solution with boundary trace $\nu$.
Furthermore, the solution is unique in the class
of \gsmod solutions. \es\medskip

\par For $q_c\leq q\leq 2$,  results similar to those stated in the last two theorems,
were obtained by Dynkin and Kuznetsov \cite{DK3} and Kuznetsov
\cite{Ku}, based on their definition of fine trace.
  However, by their results, the prescribed trace is attained
  only up to equivalence,  i.e., up to
a set of capacity zero. According to the present results, the
solution attains precisely the prescribed trace and this applies
to all $q\geq q_c$. The relation between the Dynkin-Kuznetsov
definition and the
definition presented here, is not yet clear.

\section{Main ideas of proofs} We need some additional notation.
Let $F\sbs \bdw$ be a $\GTT_q$-closed set  and let $U_F$ denote
the maximal solution vanishing on $\bdw\sms F$. Then $\inf(u,U_F)$
is a supersolution of \rqq and the largest solution dominated by
it is denoted by $[u]_F$.\1
{\em On the proof of \rth{dichotomy}.} One of the essential
features of boundary trace is its local nature. This is used in the present proof,
through the following lemmas.\medskip

\blemma{open1} Let $u\in \CU$ and
let $\set{\gb_n}$ be a \seq converging to zero such that
$w=\lim_{n\tin}u_{\gb_n}^Q$ exists,(see Notation 1.1).
 Then $[u]_F\leq w\leq [u]_{\tl Q}$, for every $\GTT_q$-closed subset of $Q$. \es\medskip
 
\blemma{open2} Let $u\in \CU$. Suppose
that there exists a $\GTT_q$-open set $Q\sbs\bdw$ and a \seq
$\set{\gb_n}$ converging to zero such that
$\sup_n\int_{\Gs_{\gb_n}(Q)}udS<\infty.$
\par Then, for any $\GTT_q$-closed set $F\sbsq Q$, $[u]_F$ is a
moderate solution.
 If $D$ is a $\GTT_q$-open set such that $\tl D\sbsq Q$, there  exists a bounded Borel
measure $\mu_D$ on $\bdw$ such that $\mu_D(\bdw\sms\tl D)=0$ and
\begin{equation}\label{regularconv}
u(\gb,\cdot)\chr{D}\rightharpoonup \mu_D \q\text{weakly relative to $C(\bdw)$ as $\gb\to 0$}.
\end{equation}
\es\medskip

\par  The dichotomy stated in the theorem is derived by combining these lemmas.\2
{\em On the proof of \rth{int-B}.} \rth{dichotomy} and
\rlemma{open2} imply that for every point $\gs\in \CR(u)$ there
exists a $\GTT_q$-open \ngh $D$ of $\gs$ and a a bounded Borel
measure $\mu_D$ on $\bdw$ such that \req{regularconv} holds.
It is not difficult to verify that, if $\gs,\gs'\in \CR(u)$ and
$D, D'$ are $\GTT_q$-open neighborhoods related to these
points as above then, the corresponding measures
$\mu\ind{D}$ and $\mu\ind{D'}$ are compatible: if $E$ is a Borel
subset of $D\cap D'$ then $\mu\ind{D}(E)=\mu\ind{D'}(E)$. The
existence of a measure $\mu^0$, which vanishes outside $\CR(u)$
and satisfies assertion (i) is easily deduced from these facts.
\par With $\gs$ and $D$ as above, if $Q$ is a $\GTT_q$-open set such that $\tl Q\sbsq D$ then $[u]_Q$ is moderate and
$\mu^0\chr{\tl Q}=\mu\ind{D}\chr{\tl Q}=\tr[u]_Q$. This implies
that $\mu^0$ is  \abc relative to $\Capq$ and outer regular. On the other
hand, if $A$ is a $\GTT_q$-open set such that $\mu^0(A)<\infty$, it
follows that on every $\GTT_q$-open set $Q$ such that $\Capq(\tl Q\sms A)=0$,
\req{intlim} holds \wrto $\mu^0$. Thus $\mu^0$ satisfies
(ii) in $\CR(u)$. The measure $\mu$  defined by
\begin{equation}\label{regularization}
  \mu(E)=\inf\set{\mu_0(D):\; \forall D\text{ \qop, } E\sbs D}
\end{equation}
satisfies (ii) on the whole boundary. It is called the $\GTT_q$-{\em regularization} of $\mu_0$.\2
{\em On the proof of \rth{exist}.} The main ingredients in this
proof are provided by \rth{int-B} (ii) and the following lemmas.\medskip

\blemma{boundary1} Let $F\sbs\bdw$ be a \qcl set. Then
$\CS(U_F)=b_q(F)$ where $b_q(F)$ denotes the set of $\Capq$-thick
points of $F$. \es \medskip

\noindent {\it Proof}  Let $\gx$ be a point on $\bdw$ such that $F$
is $\Capq$-thin at $\gx$. Then there exists a $\GTT_q$-open \ngh
$Q$ such that $\Capq(\tl Q\cap F)=0$ and \consy
$[U_F]_{Q}=U_{F\cap \tl Q}=0$. Therefore $\gx\in \CR(U_F)$.
\par Conversely, assume that $\gs\in F\cap\CR(U_F)$ and let $Q$ be a $\GTT_q$-open \ngh as in
\rth{int-B}.
 Let $D$ be a \qop \ngh
 of $\gs$ such that $\tl D\sbsq Q$. Then  $[u]_D$ is moderate and \consy
 $D\sbs \CR(u)$. In turn this implies that $\Capq(F\cap D)=0$ and \consy $F$ is $q$-thin at $\gx$.
\medskip 

\blemma{boundary2} Let $u\in \CU(\Gw)$ and let $\nu:=\tr
u$. Put
\begin{equation}\label{bnd1}
 \CS_0(u):=\set{\gx\in\bdw:\, \nu(Q\sms\CS(u))=\infty \forevery Q: \; \gx\in Q,\; Q\text{ \qop}}.
\end{equation}
Then $\CS(u)=\CS_0(u)\cup b_q(\CS(u))$. \es \medskip

The fact that
$\CS(u)\supset\CS_0(u)\cup b_q(\CS(u))$ is straightforward. The
opposite inclusion  depends on the fact that if
$\gx\in\bdw\sms b_q(\CS(u))$, there exists a $\GTT_q$-open
\ngh $Q$ such that $\Capq(\tl Q\cap \CS(u))=0$.
Hence $[u]_Q=[u]_{Q\sms \CS(u)}$. If $\gx\in
\CS(u)$ then \req{int-sing} holds. Therefore, using \req{PT},
 it is easy to show that $\nu(Q\sms\CS(u))=\infty$.\medskip
 
\blemma{gsmodtrace1}Let $u\in \CU(\Gw)$ be a \gsmod solution and 
$\set{u_n}$ an increasing sequence of moderate
solutions such that $u_n\uar u$. If $w$ is a moderate solution
dominated by $u$ then $\tr w\leq \lim \tr u_n$. \es \medskip

Put $\tau:=\tr
w$ and $\mu_0:=\lim \tr u_n$. It is sufficient to show that
$\tau(K)\leq \mu_0(K)$ for every compact set $K\sbs \bdw$ such that
$\mu_0(K)<\infty$ and $\Capq(K)>0$. It can be shown that, under
these assumptions:
 (a) There exists a $\GTT_q$-open set $Q$ such that $\Capq(K\sms Q)=0$ and $\mu_0(Q)<\infty$.
(b) If $F$ is a $\GTT_q$-closed subset of $Q$ then $\nu(F)\leq
\mu_0(Q)$.
\par These facts imply in a straightforward
manner that $\tau(K)\leq \mu_0(K)$.\2
\indent Combining  \rlemma{boundary1}, \rlemma{boundary2} and
\rth{int-B} (ii) one can verify that the existence of a solution
$u\in \CU(\Gw)$ such that $\nu=\tr u$ implies that $\nu$ is
$q$-perfect.

On the other hand, we observe that the solution $U$ defined by
\req{precisely} is \gsmod. This is based on the fact that $U_F$ is
\gsmod, which (for $F$ compact) was  established by the authors in \cite{MV4}
and it
remains valid when $F$ is $\GTT_q$-closed. Using this fact
and \rlemma{gsmodtrace1} it can be shown that, if $\nu$
is $q$-perfect then $\tr U=\nu$.\2
{\em On the proof of \rth{unique}.} The uniqueness result is
based on,\medskip

 \blemma{gsmodtrace2} Let $u\in
\CU(\Gw)$ be a \gsmod solution and $\set{u_n}$ an increasing sequence of moderate solutions such that $u_n\uar
u$. Put $\nu_0=\lim\tr u_n$ and let $\nu$ be the $\GTT_q$-regularization of
$\nu_0$. Then $\nu$ is the
precise boundary trace of $u$. In particular it is independent of
the choice of the \seq $\set{u_n}$. \es\medskip	

Let $\set{u_n}$ and $\set{w_n}$ be two increasing  {\seq}s  of
moderate solutions converging to $u$. Then \rlemma{gsmodtrace1}
implies that $\lim \tr u_n=\lim \tr w_n$. Thus $\nu$ is
independent of the choice of the \seq and it is not difficult to
verify that $\tr u=\nu$.

Clearly, \rlemma{gsmodtrace2} implies that, if $u,v$ are \gsmod
solutions and $u\leq v$ then $\tr u\leq \tr v$. Hence the
uniqueness result. In addition, if $u\in \CU(\Gw)$ and $\nu=\tr u$
then the solution $v$ defined as in \rth{exist} is uniquely
determined by $\nu$ and therefore by $u$. It can be sown that
$u\ominus v$ (= smallest solution dominating $u-v$) vanishes on $\CR(u)$. Therefore $u\ominus v\leq U_F$
where $F=\CS(u)$. Hence $U=v\oplus U_F$ is the maximal solution
with trace $\nu$.
\vskip 2mm
\noindent{\bf Acknowledgment.} Both authors were partially sponsored by an EC grant through
the RTN Program ÒFront-SingularitiesÓ, HPRN-CT-2002-00274 and by the French-Israeli cooperation program through grant
No. 3-1352. The first author (MM) also wishes to acknowledge the support of the Israeli Science Foundation through grant
No. 145-05.


\begin{thebibliography}{99}

\bibitem{AH} Adams D. R. and Hedberg L. I., Function spaces and potential theory,
Grundlehren  Math. Wissen. {\bf 314}, Springer (1996).
\bibitem{Dbook1} Dynkin E. B. {\em Diffusions, Superdiffusions and Partial Differential Equations},
American Math. Soc., Providence, Rhode Island, Colloquium
Publications {\bf 50}, 2002.
\bibitem{Dbook2}  Dynkin E. B. {\em Superdiffusions and Positive Solutions of Nonlinear
Partial Differential Equations}, American Math. Soc., Providence,
Rhode Island, Colloquium Publications {\bf 34}, 2004.
\bibitem {DK1} Dynkin E. B. and Kuznetsov S. E. {\em Superdiffusions
and removable singularities for quasilinear partial differential
equations}, Comm. Pure Appl. Math. {\bf 49}, 125-176 (1996).
\bibitem{DK3} Dynkin E. B. and Kuznetsov S. E. {\em Fine topology and fine trace on the boundary associated with
a class of quasilinear differential equations}, Comm. Pure Appl.
Math. {\bf 51}, 897-936 (1998).
\bibitem{Ku} Kuznetsov S. E. {\em \gsmod solutions of $Lu=u^\ga$ and fine trace on the boundary}, C.R.
Acad. Sc. Serie I {\bf 326}, 1189-1194 (1998).
\bibitem{LG} Legall J. F., {\em The Brownian snake and solutions of
$\Gd u=u^{2}$ in a domain}, Probab. Th. Rel. Fields {\bf 102},
393-432 (1995).
\bibitem{LGbook} Legall J. F., {\em Spatial branching processes, random snakes and partial differential
equations}, Birkh\"{a}user, Basel/Boston/Berlin, 1999.
\bibitem{MV1} Marcus M. and V\'{e}ron L., {\em The boundary trace of positive
solutions of semilinear elliptic equations: the subcritical case},
Arch. Rat. Mech. Anal. {\bf 144}, 201-231 (1998).
\bibitem{MV2} Marcus M. and V\'{e}ron L., {\em The boundary trace of positive
solutions of semilinear elliptic equations: the supercritical
case}, J. Math. Pures Appl. {\bf 77}, 481-524 (1998).
\bibitem{MV3} Marcus M. and V\'{e}ron L., {\em Removable
singularities and boundary trace}, J. Math. Pures Appl. {\bf 80},
879-900 (2000).
\bibitem{MV4} Marcus M. and V\'{e}ron L., {\em Capacitary estimates of positive solutions of semilinear
elliptic equations with absorption}, J.European Math. Soc. {\bf
6}, 483-527 (2004). {\bf }(2004).
\bibitem{Ms} Mselati B., {\em Classification and probabilistic representation of the positive solutions
of a semilinear elliptic equation}, Mem. Am. Math. Soc. {\bf 168} (2004).
\end{thebibliography}
\end{document}